\documentclass[11pt]{amsart}

\usepackage[T1]{fontenc}
\usepackage[utf8]{inputenc}
\usepackage{amsmath,amssymb,amsthm,tikz,enumitem,booktabs}
\usepackage[margin=1in]{geometry}
\usepackage[hidelinks]{hyperref}
\usepackage[title]{appendix}
\usepackage[mathlines]{lineno}
\usepackage[normalem]{ulem}
\theoremstyle{plain}
\newtheorem{theorem}{Theorem}[section]
\newtheorem{proposition}[theorem]{Proposition}
\newtheorem{lemma}[theorem]{Lemma}
\newtheorem{corollary}[theorem]{Corollary}

\theoremstyle{definition}

\newtheorem{example}[theorem]{Example}

\theoremstyle{remark}
\newtheorem{remark}[theorem]{Remark}

\numberwithin{equation}{section}

\newcommand{\Nzero}{\mathbb{Z}_{\geq 0}}
\newcommand{\fall}[2]{(#1)_{#2}}
\DeclareMathOperator{\im}{im}

\title{A General Triangle Generation of the $d$-Dimensional Genocchi Numbers}

\author{Derek Hanely}
\address{Department of Mathematics,
Penn State Behrend, Erie, Pennsylvania, USA}
\email{derek.hanely@psu.edu}

\author{Edison Mu\c{c}llari}
\email{edisonmucllari@gmail.com}

\author{Vasily Zadorozhnyy}
\address{Applied Sciences Group,
Microsoft Corporation, Redmond, Washington, USA}
\email{vzadorozhnyy@microsoft.com}


\date{}
\begin{document}

%%% pure text for arxiv and conferences/jounrals where math mode is not allowed in the abstract
% We obtain a uniform recurrence for the rook numbers of d-dimensional Genocchi boards, for every d greater than or equal to 3. The recurrence expresses each row of the rook-number triangle in terms of the preceding row, with coefficients given by products of consecutive powers and finite sums involving Stirling numbers of the second kind. The proof decomposes the boundary of the board according to the coordinates equal to 1. A coordinatewise counting argument shows that the number of extensions of a fixed interior placement depends on a boundary partition only through the number of its blocks and the size of their union. We also express the coefficients as finite differences, recover the previously established recurrences in dimensions three through five, and give the six-dimensional specialization. As a complementary result, we formulate an injection model that provides a direct bijective proof of a formerly established identification of the diagonal rook numbers with generalized Genocchi numbers.
%%%

\begin{abstract}
    We establish a uniform recurrence relation for the rook numbers of generalized Genocchi boards of dimension $d\geq 3$. The recurrence for each row of the rook-number triangle is expressed in terms of the preceding row, with coefficients given by products of consecutive powers and finite sums involving Stirling numbers of the second kind. The proof decomposes the boundary of the board into blocks that correspond to the coordinates equal to 1. A coordinatewise counting argument then shows that the number of extensions of a fixed interior placement depends on a boundary partition only through the number of its blocks and the size of their union. We also express the coefficients as finite differences, recovering the previously established recurrences in dimensions three through five, and give the six-dimensional specialization. As a complementary result, we formulate an injection model that provides a direct bijective proof of a formerly established identification of the diagonal rook numbers with generalized Genocchi numbers.
\end{abstract}

\maketitle

%----------------------------------------------------------------------
\section{Introduction}\label{sec:introduction}
%----------------------------------------------------------------------
Rook theory studies the placement of nonattacking rooks on various types of boards and relates these placements to the combinatorial-enumerative problems of permutations and special number sequences. A board in higher dimensions is a finite set of integer tuples such that two rooks on a board attack each other if their cells have at least one coordinate in common~\cite{10.2140/involve.2013.6.35}. \emph{Genocchi boards} are a special family of boards in higher dimensions formed by taking the complement of its corresponding triangular board. For a comprehensive overview of the background and definitions used in this paper, we refer the reader to~\cite{10.2140/involve.2013.6.35, Alayont2012ROOKNI}.

For $n\in\Nzero$, write $[n]=\{1,\ldots,n\}$ with $[0]=\varnothing$. For $d\geq 3$, the \emph{$d$-dimensional Genocchi board of size $n$} is
\begin{equation}\label{eq:board}
    \Gamma_n^{(d)} = \left\{
        (x_1,\ldots,x_d)\in[n]^d:
        \min_{1\leq a\leq d-1}x_a\leq x_d
    \right\}.
\end{equation}
In particular, $\Gamma_0^{(d)}=\varnothing$. A nonattacking placement is a set of cells with distinct entries in each coordinate. We denote the number of placements of $k$ nonattacking rooks on $\Gamma_n^{(d)}$ as $r_k(\Gamma_n^{(d)})$ and use the conventions
\begin{equation*}
    r_0(\Gamma_n^{(d)})=1, \qquad
    r_k(\Gamma_n^{(d)})=0 \quad\text{if } k<0 \text{ or } k>n.
\end{equation*}

A recurrence for generating the triangle of rook-numbers for three-dimensional Genocchi boards was obtained by Alayont and Krzywonos~\cite{Alayont_Krzywonos_2016}. Later, Alayont, Loewen, and Zadorozhnyy~\cite{Alayont_Loewen_Zadorozhnyy_2025} established corresponding recurrences in dimensions four and five. In each case, the recurrence expresses the rook numbers of a board of size $n$ in terms of those of size $n-1$. While the coefficients in these recurrences become increasingly complex as the dimension increases, the common structure of the recurrences suggests that there is a single counting mechanism that underlies them.

The main result of the paper is a uniform row-to-row recurrence for the rook-number triangle with explicit coefficients for every $d\geq 3$. In particular, we unify all previously known recurrences for dimensions three through five via a single, straightforward boundary-counting argument.

We use the standard notation to denote the Stirling number of the second kind, $S(a,b)$, i.e., the number of ways to partition a set of $a$ elements into $b$ nonempty subsets. For $b<0$ and $b>a$, we set $S(a,b)=0$, and we use the standard convention for the empty product by setting $S(0,0)=1$.

\begin{theorem}[Uniform rook-number recurrence]\label{thm:main}
    Let $d\geq 3$, $n\geq 1$, and $0\leq k\leq n$, and put $m=n-k+1$. Then,
    \begin{equation}\label{eq:main-recurrence}
        r_k(\Gamma_n^{(d)}) = \sum_{j=0}^{\min(k,d-1)}C_j^{(d)}(m)\,r_{k-j}(\Gamma_{n-1}^{(d)}),
    \end{equation}
    where $C_0^{(d)}(m)=1$ and, for $1\leq j\leq d-1$,
    \begin{equation}\label{eq:main-coefficients}
        C_j^{(d)}(m) = (m+j-1) \prod_{a=0}^{j-2}(m+a)^d \sum_{\sigma=j}^{d-1} \binom{d-1}{\sigma} S(\sigma,j)(m-1)^{d-1-\sigma}.
    \end{equation}
\end{theorem}

The proof separates the contribution of a boundary placement into two parts. First, we will count the number of possible sets of coordinates that are equal to~1 for a fixed boundary placement. Then, we calculate the number of possible assignments of the remaining coordinates. For a fixed $j$, this count will depend on the total number of coordinates equal to~1 and not on how they are distributed among the rooks; this justifies why the sum in~\eqref{eq:main-coefficients} is a sum of Stirling numbers of the second kind.

The diagonal sequence $r_n(\Gamma_n^{(d)})$ has already been identified with generalized Genocchi numbers by Alayont, Moger-Reischer, and Swift~\cite{Alayont2012ROOKNI}. Their work also gives a combinatorial interpretation of the generalized Genocchi numbers as $(d-1)$-tuples of permutations. As a complementary result, in Section~\ref{sec:gandhi}, we give a direct and intuitive bijective proof of this established identification. We also enlarge the permutation model to allow for injections into larger sets, which gives a simple interpretation of the subtraction term in the generalized Gandhi recurrence by deleting a chosen unused value from each codomain.

The paper is organized as follows. We first identify the boundary-counting mechanism underlying the uniform recurrence in Section~\ref{sec:proof}. Then, we explore the algebraic properties of the coefficients in Section~\ref{sec:coefficients}, and relate the general formula to some explicit examples in low dimensions in Section~\ref{sec:examples}. We move on from the full rook-number triangle in Section~\ref{sec:gandhi} to provide a bijective explanation of the diagonal entries. Finally, we discuss some further questions in Section~\ref{sec:conclusion}, and list some more numerical examples for various dimensions in Appendix~\ref{appendix_coefficients}.

%----------------------------------------------------------------------
\section{Boundary decomposition and proof of the recurrence}\label{sec:proof}
%----------------------------------------------------------------------
The recurrence in Theorem~\ref{thm:main} is derived by first partitioning the original board into a smaller translated board and its remaining boundary, then counting the number of placements on the smaller board that extend across that boundary.
\subsection{The translated interior}

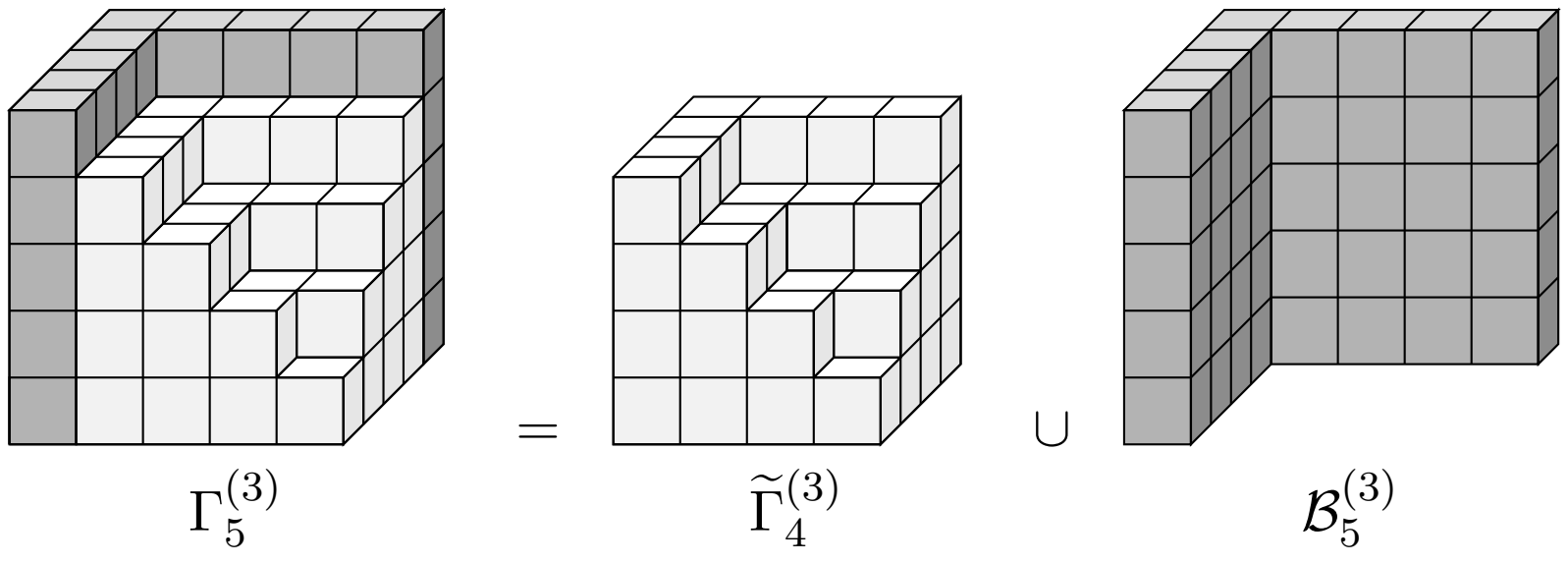
\begin{figure}[t]
    \centering
    \begingroup
    \setlength{\tabcolsep}{5pt}
    \renewcommand{\arraystretch}{1.15}

    % Overlapping hyperplanes.
    \begin{tabular}{@{}c c c c c@{}}
        \begin{tikzpicture}[scale=0.45]
% 0 Rooks outside, all inside n-1 size board
% \draw[fill=black!30]  (0,0)--(0,5)--(1.5,6.5)--(6.5,6.5)--(6.5,1.5)--(6.2,1.2)--(6.2,5.2)--(2.2,5.2)--(1,4)--(1,0)--cycle;

\draw[fill=black!30]  
    (0,0)--(0,5)--(1,5)--(1,0)--cycle
    (2.2,5.2)--(2.2,6.2)--(6.2,6.2)--(6.2,5.2)--cycle
;
\draw[fill=black!15] (0,5)--(1.5,6.5)--(6.5,6.5)--(6.2,6.2)--(2.2,6.2)--(1,5)--cycle;
\draw[fill=black!45] 
    (6.5,6.5)--(6.5,1.5)--(6.2,1.2)--(6.2,6.2)--cycle
    (1,4)--(1,5)--(2.2,6.2)--(2.2,5.2)--cycle
;

\draw[fill=black!10]  
    (5,0)--(6.2,1.2)--(6.2,5.2)--(5.9,4.9)--(5.9,3.9)--(5.6,3.6)--(5.6,2.6)--(5.3,2.3)--(5.3,1.3)--(5,1)--cycle
    (4,1)--(4,2)--(4.3,2.3)--(4.3,1.3)--cycle
    (3,2)--(3,3)--(3.6,3.6)--(3.6,2.6)--cycle
    (2,3)--(2,4)--(2.9,4.9)--(2.9,3.9)--cycle
;
\draw[fill=black!5]
    (1,0)--(1,4)--(2,4)--(2,3)--(3,3)--(3,2)--(4,2)--(4,1)--(5,1)--(5,0)--cycle
    (4.3,1.3)--(4.3,2.3)--(5.3,2.3)--(5.3,1.3)--cycle
    (3.6,2.6)--(3.6,3.6)--(5.6,3.6)--(5.6,2.6)--cycle
    (2.9,3.9)--(2.9,4.9)--(5.9,4.9)--(5.9,3.9)--cycle
;

\draw (0,0)--(5,0)--(5,1)--(0,1)--(0,0);
\draw (0,2)--(4,2)--(4,0);
\draw (0,3)--(3,3)--(3,0);
\draw (0,4)--(2,4)--(2,0);
\draw (0,5)--(1,5)--(1,0);
\draw (0,5)--(0,0);
\draw (5,1)--(5.3,1.3)--(5.3,0.3)--(5,0);
\draw (5.3,1.3)--(5.3,2.3)--(4.3,2.3)--(4,2);
\draw (4,1)--(4.3,1.3)--(5.3,1.3);
\draw (4.3,1.3)--(4.3,2.3);
\draw (3,2)--(3.6,2.6)--(5.6,2.6);
\draw (3,3)--(3.6,3.6)--(5.6,3.6);
\draw (2,4)--(2.9,4.9)--(5.9,4.9);
\draw (1,5)--(2.2,6.2)--(6.2,6.2);
\draw (5.3,2.3)--(5.6,2.6)--(5.6,3.6)--(5.9,3.9)--(5.9,4.9)--(6.2,5.2)
--(6.2,6.2)--(6.5,6.5);
\draw (0,5)--(1.5,6.5)--(6.5,6.5);
\draw (1,4)--(2.2,5.2)--(6.2,5.2)--(6.5,5.5);
\draw (2,3)--(2.9,3.9)--(5.9,3.9);
\draw (0.3,5.3)--(1.3,5.3)--(1.3,4.3)--(2.3,4.3)--(2.3,3.3)--(3.3,3.3)
--(3.3,2.3)--(4.3,2.3)--(4.3,1.3)--(5.3,1.3);
\draw (0.6,5.6)--(1.6,5.6)--(1.6,4.6)--(2.6,4.6)--(2.6,3.6)--(3.6,3.6)
--(3.6,2.6);
\draw (0.9,5.9)--(1.9,5.9)--(1.9,4.9)--(2.9,4.9)--(2.9,3.9);
\draw (1.2,6.2)--(2.2,6.2)--(2.2,5.2);
\draw (4.3,2.3)--(4.6,2.6)--(4.6,3.6)--(4.9,3.9)--(4.9,4.9)--(5.2,5.2)--
(5.2,6.2)--(5.5,6.5);
\draw (3.6,3.6)--(3.9,3.9)--(3.9,4.9)--(4.2,5.2)--(4.2,6.2)--(4.5,6.5);
\draw (2.9,4.9)--(3.2,5.2)--(3.2,6.2)--(3.5,6.5);
\draw (2.2,6.2)--(2.5,6.5);
\draw (6.5,6.5)--(6.5,1.5)--(5.3,0.3);
\draw (5.3,1.3)--(6.5,2.5);
\draw (5.6,2.6)--(6.5,3.5);
\draw (5.9,3.9)--(6.5,4.5);
\draw (6.2,4.2)--(6.5,4.5);
\draw (6.2,5.2)--(6.2,1.2);
\draw (5.3,1.3)--(5.3,0.3);
\draw (5.6,2.6)--(5.6,0.6);
\draw (5.9,3.9)--(5.9,.9);
\end{tikzpicture}
        & $=$
        & \begin{tikzpicture}[scale=0.45]

\draw[fill=black!10] 
    (4,0)--(5.2,1.2)--(5.2,5.2)--(4.9,4.9)--(4.9,3.9)--(4.6,3.6)--(4.6,2.6)--(4.3,2.3)--(4.3,1.3)--(4,1)--cycle
    (3,1)--(3,2)--(3.3,2.3)--(3.3,1.3)--cycle
    (2,2)--(2,3)--(2.6,3.6)--(2.6,2.6)--cycle
    (1,3)--(1,4)--(1.9,4.9)--(1.9,3.9)--cycle
;

\draw[fill=black!5]
    (0,0)--(0,4)--(1,4)--(1,3)--(2,3)--(2,2)--(3,2)--(3,1)--(4,1)--(4,0)--cycle
    (3.3,1.3)--(3.3,2.3)--(4.3,2.3)--(4.3,1.3)--cycle
    (2.6,2.6)--(2.6,3.6)--(4.6,3.6)--(4.6,2.6)--cycle
    (1.9,3.9)--(1.9,4.9)--(4.9,4.9)--(4.9,3.9)--cycle
;

\draw (0,0)--(4,0)--(5.2,1.2)--(5.2,5.2)--(1.2,5.2)--(0,4)--cycle;
\draw (1,0)--(1,4)--(1.9,4.9)--(4.9,4.9)--(4.9,.9);
\draw (1,3)--(1.9,3.9)--(4.9,3.9);
\draw (2,0)--(2,3)--(2.6,3.6)--(4.6,3.6)--(4.6,.6);
\draw (2,2)--(2.6,2.6)--(4.6,2.6);
\draw (3,0)--(3,2)--(3.3,2.3)--(4.3,2.3)--(4.3,.3);
\draw (0,1)--(4,1)--(5.2,2.2);
\draw (0,2)--(3,2);
\draw (0,3)--(2,3);
\draw (0,4)--(1,4);
\draw (1.9,4.9)--(2.2,5.2);
\draw (3,1)--(3.3,1.3)--(4.3,1.3);
\draw (0.3,4.3)--(1.3,4.3)--(1.3,3.3)--(2.3,3.3)--(2.3,2.3)--(3.3,2.3);
\draw (0.6,4.6)--(1.6,4.6)--(1.6,3.6)--(2.6,3.6)--(2.6,2.6);
\draw (0.9,4.9)--(1.9,4.9)--(1.9,3.9);
\draw (3.2,5.2)--(2.9,4.9)--(2.9,3.9)--(2.6,3.6);
\draw (4.2,5.2)--(3.9,4.9)--(3.9,3.9)--(3.6,3.6)--(3.6,2.6)--(3.3,2.3)--(3.3,1.3);
\draw (4,0)--(4,1);
\draw (4.3,2.3)--(5.2,3.2);
\draw (4.6,3.6)--(5.2,4.2);
\draw (4.9,4.9)--(5.2,5.2);
\end{tikzpicture}
        & $\cup$
        & \begin{tikzpicture}[scale=0.45]
% 0 Rooks outside, all inside n-1 size board
% \draw[fill=black!30]  (0,0)--(0,5)--(1.5,6.5)--(6.5,6.5)--(6.5,1.5)--(6.2,1.2)--(2.2,1.2)--(1,0)--cycle;

\draw[fill=black!30]  (0,0)--(0,5)--(1,5)--(1,0)--cycle;
\draw[fill=black!30]  (2.2,1.2)--(2.2,6.2)--(6.2,6.2)--(6.2,1.2)--cycle;
\draw[fill=black!15] (0,5)--(1.5,6.5)--(6.5,6.5)--(6.2,6.2)--(2.2,6.2)--(1,5)--cycle;
\draw[fill=black!45] (6.5,6.5)--(6.5,1.5)--(6.2,1.2)--(6.2,6.2)--cycle;
\draw[fill=black!45] (1,0)--(1,5)--(2.2,6.2)--(2.2,1.2)--cycle;

\draw (1,0)--(1,5)--(1,5)--(2.2,6.2)--(6.2,6.2)--(6.5,6.5);
\draw (1.3,0.3)--(1.3,5.3)--(0.3,5.3);
\draw (1.6,0.6)--(1.6,5.6)--(0.6,5.6);
\draw (1.9,0.9)--(1.9,5.9)--(0.9,5.9);
\draw (2.2,1.2)--(2.2,6.2)--(1.2,6.2);

\draw (3.2,1.2)--(3.2,6.2)--(3.5,6.5);
\draw (4.2,1.2)--(4.2,6.2)--(4.5,6.5);
\draw (5.2,1.2)--(5.2,6.2)--(5.5,6.5);
\draw (6.2,1.2)--(6.2,6.2);

\draw (0,1)--(1,1)--(2.2,2.2)--(6.2,2.2)--(6.5,2.5);
\draw (0,2)--(1,2)--(2.2,3.2)--(6.2,3.2)--(6.5,3.5);
\draw (0,3)--(1,3)--(2.2,4.2)--(6.2,4.2)--(6.5,4.5);
\draw (0,4)--(1,4)--(2.2,5.2)--(6.2,5.2)--(6.5,5.5);

\draw (2.2,6.2)--(2.5,6.5);
\draw (0,5)--(1,5);
\end{tikzpicture}
        \\
        $\Gamma_5^{(3)}$
        &&
        $\widetilde{\Gamma}_4^{(3)}$
        &&
        $\mathcal{B}_5^{(3)}$
    \end{tabular}

    \endgroup

    \caption{Decomposition of the Genocchi board $\Gamma_5^{(3)}$ into its interior - $\widetilde{\Gamma}_4^{(3)}$ and boundary $\mathcal{B}_5^{(3)}=\Gamma_5^{(3)} \setminus\widetilde{\Gamma}_4^{(3)}$ parts, respectively, using the visualization from~\cite{Alayont_Krzywonos_2016}.}
    \label{fig:board-tranlation-boundary}
\end{figure}

To compare successive board sizes, we use the translated copy
\begin{equation*}%\label{eq:interior}
    \widetilde{\Gamma}_{n-1}^{(d)} = \left\{
        (x_1,\ldots,x_d)\in\Gamma_n^{(d)}: x_a\geq 2 \text{ for every }a \in [d]
    \right\}.
\end{equation*}
Translation by $(-1,\ldots,-1)$ gives a bijection from $\widetilde{\Gamma}_{n-1}^{(d)}$ to $\Gamma_{n-1}^{(d)}$ which also preserves nonattacking rook placements. Thus,
\begin{equation*}
    r_\ell(\widetilde{\Gamma}_{n-1}^{(d)}) = r_\ell(\Gamma_{n-1}^{(d)}).
\end{equation*}

Now, define
\begin{equation*}
    \mathcal{H}_a = \{x\in\Gamma_n^{(d)}:x_a=1\} \text{ for } a\in [d-1].
\end{equation*}
It follows that the boundary added at size $n$ is precisely
\begin{equation}\label{eq:boundary}
    \mathcal{B}_n^{(d)} := \Gamma_n^{(d)} \setminus\widetilde{\Gamma}_{n-1}^{(d)} = \bigcup_{a\in [d-1]}\mathcal{H}_a.
\end{equation}
Note that every cell outside the translated interior has at least one coordinate equal to~1. If this coordinate is $x_d$, the board condition implies that there must exist some $a\in [d-1]$ such that $x_a=1$ as well. Conversely, every cell in an $\mathcal{H}_a$ lies entirely outside the interior. An example of the decomposition of the board into both interior and boundary components for $d=3$ and $n=5$ is given in Figure~\ref{fig:board-tranlation-boundary}.

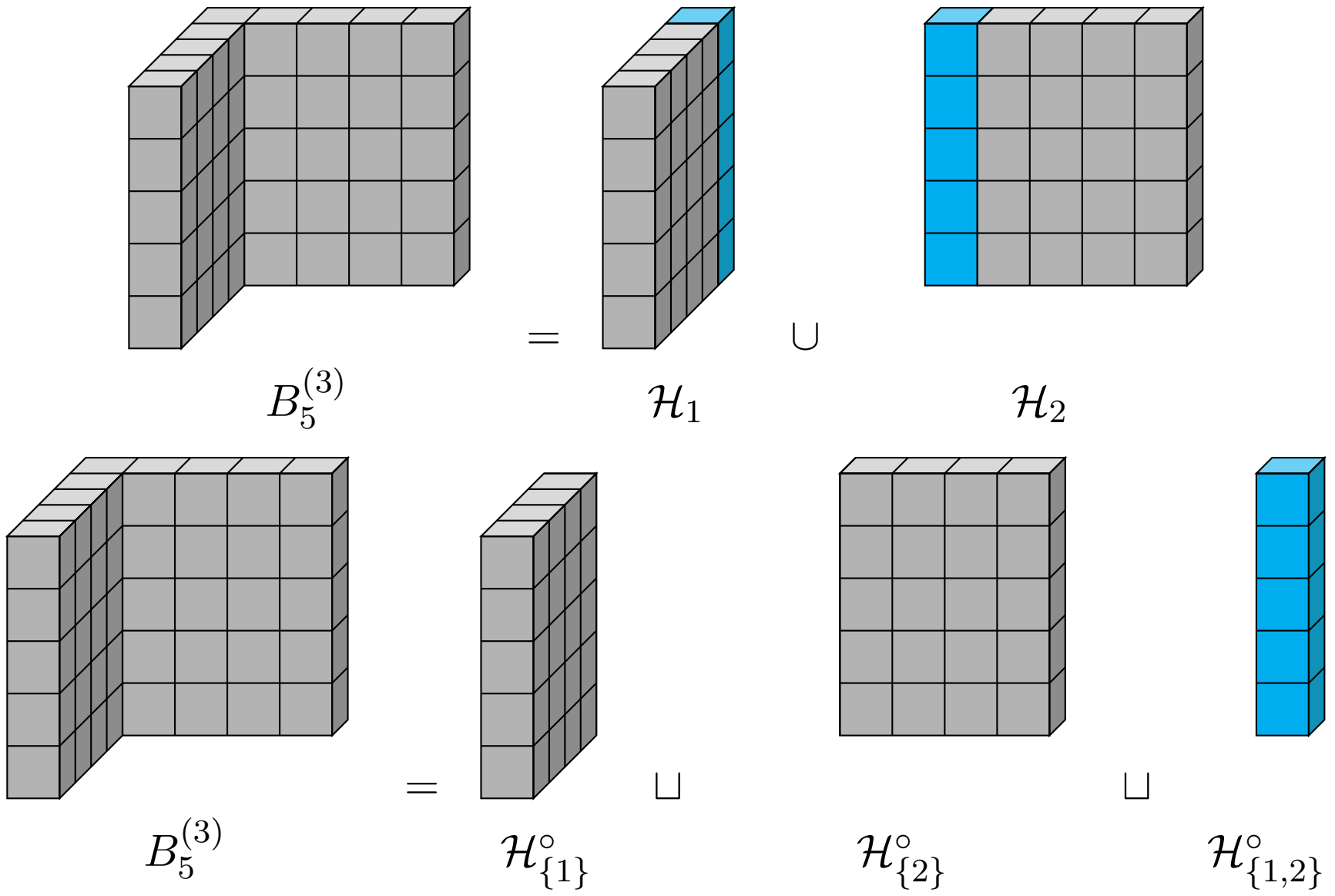
\begin{figure}[t]
    \centering
    \begingroup
    \setlength{\tabcolsep}{5pt}
    \renewcommand{\arraystretch}{1.15}

    % Overlapping hyperplanes.
    \begin{tabular}{@{}c c c c c@{}}
        \begin{tikzpicture}[scale=0.45]
% 0 Rooks outside, all inside n-1 size board
% \draw[fill=black!30]  (0,0)--(0,5)--(1.5,6.5)--(6.5,6.5)--(6.5,1.5)--(6.2,1.2)--(2.2,1.2)--(1,0)--cycle;

\draw[fill=black!30]  (0,0)--(0,5)--(1,5)--(1,0)--cycle;
\draw[fill=black!30]  (2.2,1.2)--(2.2,6.2)--(6.2,6.2)--(6.2,1.2)--cycle;
\draw[fill=black!15] (0,5)--(1.5,6.5)--(6.5,6.5)--(6.2,6.2)--(2.2,6.2)--(1,5)--cycle;
\draw[fill=black!45] (6.5,6.5)--(6.5,1.5)--(6.2,1.2)--(6.2,6.2)--cycle;
\draw[fill=black!45] (1,0)--(1,5)--(2.2,6.2)--(2.2,1.2)--cycle;

\draw (1,0)--(1,5)--(1,5)--(2.2,6.2)--(6.2,6.2)--(6.5,6.5);
\draw (1.3,0.3)--(1.3,5.3)--(0.3,5.3);
\draw (1.6,0.6)--(1.6,5.6)--(0.6,5.6);
\draw (1.9,0.9)--(1.9,5.9)--(0.9,5.9);
\draw (2.2,1.2)--(2.2,6.2)--(1.2,6.2);

\draw (3.2,1.2)--(3.2,6.2)--(3.5,6.5);
\draw (4.2,1.2)--(4.2,6.2)--(4.5,6.5);
\draw (5.2,1.2)--(5.2,6.2)--(5.5,6.5);
\draw (6.2,1.2)--(6.2,6.2);

\draw (0,1)--(1,1)--(2.2,2.2)--(6.2,2.2)--(6.5,2.5);
\draw (0,2)--(1,2)--(2.2,3.2)--(6.2,3.2)--(6.5,3.5);
\draw (0,3)--(1,3)--(2.2,4.2)--(6.2,4.2)--(6.5,4.5);
\draw (0,4)--(1,4)--(2.2,5.2)--(6.2,5.2)--(6.5,5.5);

\draw (2.2,6.2)--(2.5,6.5);
\draw (0,5)--(1,5);
\end{tikzpicture}
        & $=$
        & \begin{tikzpicture}[scale=0.45]
% 0 Rooks outside, all inside n-1 size board
% \draw[fill=black!30]  (0,0)--(0,5)--(1.5,6.5)--(2.5,6.5)--(2.5,1.5)--(1,0)--cycle;
\draw[fill=black!30]  (0,0)--(0,5)--(1,5)--(1,0)--cycle;
\draw[fill=black!45] (1,0)--(1,5)--(2.2,6.2)--(2.2,1.2)--cycle;
\draw[fill=black!15] (0,5)--(1.2,6.2)--(2.2,6.2)--(1,5)--cycle;

% \draw[fill=cyan] (1.2,6.2)--(1.5,6.5)--(2.5,6.5)--(2.5,1.5)--(2.2,1.2)--(2.2,6.2)--cycle;

\draw[fill=cyan!50] (1.2,6.2)--(1.5,6.5)--(2.5,6.5)--(2.2,6.2)--cycle;
\draw[fill=cyan!70!black] (2.2,6.2)--(2.5,6.5)--(2.5,1.5)--(2.2,1.2)--cycle;

\draw (1,0)--(1,5)--(2.2,6.2);
\draw (1.3,0.3)--(1.3,5.3)--(0.3,5.3);
\draw (1.6,0.6)--(1.6,5.6)--(0.6,5.6);
\draw (1.9,0.9)--(1.9,5.9)--(0.9,5.9);
\draw (2.2,1.2)--(2.2,6.2)--(1.2,6.2);

% \draw (3.2,1.2)--(3.2,6.2)--(3.5,6.5);
% \draw (4.2,1.2)--(4.2,6.2)--(4.5,6.5);
% \draw (5.2,1.2)--(5.2,6.2)--(5.5,6.5);
% \draw (6.2,1.2)--(6.2,6.2);

\draw (0,1)--(1,1)--(2.5,2.5);
\draw (0,2)--(1,2)--(2.5,3.5);
\draw (0,3)--(1,3)--(2.5,4.5);
\draw (0,4)--(1,4)--(2.5,5.5);

\draw (2.2,6.2)--(2.5,6.5);
\draw (0,5)--(1,5);
\end{tikzpicture}
        & $\cup$
        & \begin{tikzpicture}[scale=0.45]
% 0 Rooks outside, all inside n-1 size board
\draw[white] (0,0);
% \draw[fill=black!30]  (2.2,6.2)--(2.5,6.5)--(6.5,6.5)--(6.5,1.5)--(6.2,1.2)--(2.2,1.2)--cycle;
\draw[fill=black!15]  (2.2,6.2)--(2.5,6.5)--(6.5,6.5)--(6.2,6.2)--cycle;
\draw[fill=black!30]  (2.2,6.2)--(6.2,6.2)--(6.2,1.2)--(2.2,1.2)--cycle;
\draw[fill=black!45] (6.2,6.2)--(6.5,6.5)--(6.5,1.5)--(6.2,1.2)--cycle;

% \draw[fill=cyan] (1.2,6.2)--(1.5,6.5)--(2.5,6.5)--(2.2,6.2)--(2.2,1.2)--(1.2,1.2)--cycle;
\draw[fill=cyan] (1.2,1.2)--(1.2,6.2)--(2.2,6.2)--(2.2,1.2)--cycle;
\draw[fill=cyan!50] (1.2,6.2)--(1.5,6.5)--(2.5,6.5)--(2.2,6.2)--cycle;

% \draw (1,0)--(1,5)--(1,5)--(2.2,6.2)--(6.2,6.2)--(6.5,6.5);
% \draw (1.3,0.3)--(1.3,5.3)--(0.3,5.3);
% \draw (1.6,0.6)--(1.6,5.6)--(0.6,5.6);
% \draw (1.9,0.9)--(1.9,5.9)--(0.9,5.9);
\draw (2.2,1.2)--(2.2,6.2)--(2.5,6.5);

\draw (3.2,1.2)--(3.2,6.2)--(3.5,6.5);
\draw (4.2,1.2)--(4.2,6.2)--(4.5,6.5);
\draw (5.2,1.2)--(5.2,6.2)--(5.5,6.5);
\draw (6.2,1.2)--(6.2,6.2)--(6.5,6.5);

\draw (1.2,2.2)--(6.2,2.2)--(6.5,2.5);
\draw (1.2,3.2)--(6.2,3.2)--(6.5,3.5);
\draw (1.2,4.2)--(6.2,4.2)--(6.5,4.5);
\draw (1.2,5.2)--(6.2,5.2)--(6.5,5.5);
\draw (1.2,6.2)--(6.2,6.2);

\end{tikzpicture}
        \\
        $B_5^{(3)}$
        &&
        $\mathcal{H}_1$
        &&
        $\mathcal{H}_2$
    \end{tabular}

    \par\medskip

    % Pairwise disjoint exclusive regions.
    \begin{tabular}{@{}c c c c c c c@{}}
        \begin{tikzpicture}[scale=0.45]
% 0 Rooks outside, all inside n-1 size board
% \draw[fill=black!30]  (0,0)--(0,5)--(1.5,6.5)--(6.5,6.5)--(6.5,1.5)--(6.2,1.2)--(2.2,1.2)--(1,0)--cycle;

\draw[fill=black!30]  (0,0)--(0,5)--(1,5)--(1,0)--cycle;
\draw[fill=black!30]  (2.2,1.2)--(2.2,6.2)--(6.2,6.2)--(6.2,1.2)--cycle;
\draw[fill=black!15] (0,5)--(1.5,6.5)--(6.5,6.5)--(6.2,6.2)--(2.2,6.2)--(1,5)--cycle;
\draw[fill=black!45] (6.5,6.5)--(6.5,1.5)--(6.2,1.2)--(6.2,6.2)--cycle;
\draw[fill=black!45] (1,0)--(1,5)--(2.2,6.2)--(2.2,1.2)--cycle;

\draw (1,0)--(1,5)--(1,5)--(2.2,6.2)--(6.2,6.2)--(6.5,6.5);
\draw (1.3,0.3)--(1.3,5.3)--(0.3,5.3);
\draw (1.6,0.6)--(1.6,5.6)--(0.6,5.6);
\draw (1.9,0.9)--(1.9,5.9)--(0.9,5.9);
\draw (2.2,1.2)--(2.2,6.2)--(1.2,6.2);

\draw (3.2,1.2)--(3.2,6.2)--(3.5,6.5);
\draw (4.2,1.2)--(4.2,6.2)--(4.5,6.5);
\draw (5.2,1.2)--(5.2,6.2)--(5.5,6.5);
\draw (6.2,1.2)--(6.2,6.2);

\draw (0,1)--(1,1)--(2.2,2.2)--(6.2,2.2)--(6.5,2.5);
\draw (0,2)--(1,2)--(2.2,3.2)--(6.2,3.2)--(6.5,3.5);
\draw (0,3)--(1,3)--(2.2,4.2)--(6.2,4.2)--(6.5,4.5);
\draw (0,4)--(1,4)--(2.2,5.2)--(6.2,5.2)--(6.5,5.5);

\draw (2.2,6.2)--(2.5,6.5);
\draw (0,5)--(1,5);
\end{tikzpicture}
        & $=$
        & \begin{tikzpicture}[scale=0.45]
% 0 Rooks outside, all inside n-1 size board
% \draw[fill=black!30]  (0,0)--(0,5)--(1.2,6.2)--(2.2,6.2)--(2.2,1.2)--(1,0)--cycle;
\draw[fill=black!30]  (0,0)--(0,5)--(1,5)--(1,0)--cycle;
\draw[fill=black!45] (1,0)--(1,5)--(2.2,6.2)--(2.2,1.2)--cycle;
\draw[fill=black!15] (0,5)--(1.2,6.2)--(2.2,6.2)--(1,5)--cycle;

\draw (1,0)--(1,5)--(2.2,6.2);
\draw (1.3,0.3)--(1.3,5.3)--(0.3,5.3);
\draw (1.6,0.6)--(1.6,5.6)--(0.6,5.6);
\draw (1.9,0.9)--(1.9,5.9)--(0.9,5.9);
\draw (2.2,1.2)--(2.2,6.2)--(1.2,6.2);

\draw (0,1)--(1,1)--(2.2,2.2);
\draw (0,2)--(1,2)--(2.2,3.2);
\draw (0,3)--(1,3)--(2.2,4.2);
\draw (0,4)--(1,4)--(2.2,5.2);

\draw (2.2,6.2)--(2.2,6.2);
\draw (0,5)--(1,5);
\end{tikzpicture}
        & $\sqcup$
        & \begin{tikzpicture}[scale=0.45]
% 0 Rooks outside, all inside n-1 size board
\draw[white] (0,0);
% \draw[fill=black!30]  (2.2,6.2)--(2.5,6.5)--(6.5,6.5)--(6.5,1.5)--(6.2,1.2)--(2.2,1.2)--cycle;
\draw[fill=black!15]  (2.2,6.2)--(2.5,6.5)--(6.5,6.5)--(6.2,6.2)--cycle;
\draw[fill=black!30]  (2.2,6.2)--(6.2,6.2)--(6.2,1.2)--(2.2,1.2)--cycle;
\draw[fill=black!45] (6.2,6.2)--(6.5,6.5)--(6.5,1.5)--(6.2,1.2)--cycle;

\draw (2.2,1.2)--(2.2,6.2)--(2.5,6.5);

\draw (3.2,1.2)--(3.2,6.2)--(3.5,6.5);
\draw (4.2,1.2)--(4.2,6.2)--(4.5,6.5);
\draw (5.2,1.2)--(5.2,6.2)--(5.5,6.5);
\draw (6.2,1.2)--(6.2,6.2)--(6.5,6.5);

\draw (2.2,2.2)--(6.2,2.2)--(6.5,2.5);
\draw (2.2,3.2)--(6.2,3.2)--(6.5,3.5);
\draw (2.2,4.2)--(6.2,4.2)--(6.5,4.5);
\draw (2.2,5.2)--(6.2,5.2)--(6.5,5.5);
\draw (2.2,6.2)--(6.2,6.2);

\end{tikzpicture}
        & $\sqcup$
        & \begin{tikzpicture}[scale=0.45]
% 0 Rooks outside, all inside n-1 size board
\draw[white] (0,0);

% \draw[fill=cyan] (1.2,1.2)--(1.2,6.2)--(1.5,6.5)--(2.5,6.5)--(2.5,1.5)--(2.2,1.2)--cycle;
\draw[fill=cyan] (1.2,1.2)--(1.2,6.2)--(2.2,6.2)--(2.2,1.2)--cycle;
\draw[fill=cyan!50] (1.2,6.2)--(1.5,6.5)--(2.5,6.5)--(2.2,6.2)--cycle;
\draw[fill=cyan!70!black] (2.2,6.2)--(2.5,6.5)--(2.5,1.5)--(2.2,1.2)--cycle;

\draw (2.2,1.2)--(2.2,6.2)--(2.5,6.5);

\draw (1.2,2.2)--(2.2,2.2)--(2.5,2.5);
\draw (1.2,3.2)--(2.2,3.2)--(2.5,3.5);
\draw (1.2,4.2)--(2.2,4.2)--(2.5,4.5);
\draw (1.2,5.2)--(2.2,5.2)--(2.5,5.5);
\draw (1.2,6.2)--(2.2,6.2);

\end{tikzpicture}
        \\
        $B_5^{(3)}$
        &&
        $\mathcal{H}_{\{1\}}^\circ$
        &&
        $\mathcal{H}_{\{2\}}^\circ$
        &&
        $\mathcal{H}_{\{1,2\}}^\circ$
    \end{tabular}

    \endgroup

    \caption{Boundary decomposition of the three-dimensional Genocchi board of size~$5$, where $B_5^{(3)}=\Gamma_5^{(3)} \setminus\widetilde{\Gamma}_4^{(3)}$. \emph{Top:} The boundary consists of the union of two overlapping hyperplanes $\mathcal{H}_1$ and $\mathcal{H}_2$. \emph{Bottom:} By separating their intersection, the three pairwise disjoint exclusive regions $\mathcal{H}_{\{1\}}^\circ$, $\mathcal{H}_{\{2\}}^\circ$, and $\mathcal{H}_{\{1,2\}}^\circ$ are obtained. 
    The visible faces of the common intersection $\mathcal{H}_1\cap \mathcal{H}_2=\mathcal{H}_{\{1,2\}}^\circ$ are highlighted in \textcolor{cyan}{\textbf{cyan}} for each of the single hyperplanes and the intersection.}
    \label{fig:boundary-decomposition}
\end{figure}

For each $\emptyset\neq S \subseteq [d-1]$, define the \emph{exclusive region}
\begin{equation*}%\label{eq:exclusive}
    \mathcal{H}_S^\circ = \left\{x \in \Gamma_n^{(d)}: x_a=1 \text{ for }a\in S \text{ and } x_a\neq 1 \text{ for }a\in [d-1]\setminus S\right\}.
\end{equation*}
Since $S$ is nonempty, every cell in this set satisfies the board condition. Therefore, in addition to~\eqref{eq:boundary}, the boundary added at size $n$ can be represented as 
\begin{equation*}%\label{eq:exclusive-decomposition}
    \mathcal{B}_n^{(d)} = \bigsqcup_{\varnothing\neq S\subseteq [d-1]}\mathcal{H}_S^\circ.
\end{equation*}
The decomposition of the boundary into regions for $d=3$ and $n=5$ is shown in Figure~\ref{fig:boundary-decomposition}.

According to the definition of the board~\eqref{eq:board}, if any of the first $d-1$ coordinates are equal to~1, then the board condition is satisfied automatically. This means that once the exclusive region of a boundary rook is specified, there are no additional inequalities to be satisfied by the free coordinates of that rook.

\subsection{Partition types}
The exclusive regions describe where a single boundary rook can lie. To count several boundary rooks, we must determine which of these regions can occur together in a nonattacking placement. For a boundary cell $x$, let
\begin{equation*}
    S_x=\{a\in [d-1]:x_a=1\}.
\end{equation*}
Note that this set is nonempty. Moreover, if $x$ and $y$ belong to a nonattacking placement, then $S_x\cap S_y=\varnothing$. Thus, a placement of $j$ boundary rooks determines an unordered family of pairwise disjoint nonempty subsets $\mathcal{S}=\{S_1,\ldots,S_j\}$
of $[d-1]$ (in particular, $j\leq d-1$). We call such $\mathcal{S}$ its \emph{partition type}.

Write $U=S_1\cup\cdots\cup S_j$ and denote $\sigma=|U|$. Then, for fixed $\sigma$ and $j$, the number of partition types is given by
\begin{equation}\label{eq:partition-count}
    \binom{d-1}{\sigma}S(\sigma,j),
\end{equation}
where we first choose $U$ and then partition it into~$j$ nonempty blocks.

\begin{example}
    When $d=3$, the boundary has the three exclusive regions
    \begin{equation*}
        \mathcal{H}_{\{1\}}^\circ,\qquad
        \mathcal{H}_{\{2\}}^\circ,\qquad
        \mathcal{H}_{\{1,2\}}^\circ.
    \end{equation*}
    A single boundary rook can occupy any one of these regions. For two boundary rooks, the only possible partition type is $\{\{1\},\{2\}\}$.
\end{example}

\subsection{The independence lemma}

The partition count is determined by specifying the coordinates equal to~1. The remaining coordinate values can then be assigned in any way that does not create attacks. The following lemma shows that the partition count is independent of the interior placement and is determined solely by the number of blocks in the partition type and the size of the union of its blocks.

For a nonnegative integer $t$, we use the standard notation for the falling factorial:
\begin{equation*} %\label{eq:z_t}
    \fall{z}{t}=z(z-1)\cdots(z-t+1), \text{ with }\fall{z}{0}=1.
\end{equation*}

\begin{lemma}[Independence lemma]\label{lem:independence}
    Let $d\geq 3$, $n\geq 1$, $1\leq k\leq n$, and $1\leq j\leq\min(k,d-1)$. Fix a placement of $k-j$ nonattacking rooks on $\widetilde{\Gamma}_{n-1}^{(d)}$ and a partition type $\mathcal{S}=\{S_1,\ldots,S_j\}$ with $|S_1\cup\cdots\cup S_j|=\sigma$. Take $m=n-k+1$. The number of extensions with one boundary rook in each $\mathcal{H}_{S_i}^\circ$ is
    \begin{equation}\label{eq:independence}
        P_j^{(d)}(m,\sigma) = (m+j-1)\prod_{a=0}^{j-2}(m+a)^d(m-1)^{d-1-\sigma}.
    \end{equation}
    In particular, the count does not depend on the choice of the interior placement nor on the individual blocks of $\mathcal{S}$, except through $j$ and $\sigma$.
\end{lemma}

\begin{proof}
Index new rooks by their distinct blocks $S_1,\ldots,S_j$ in which they are placed. Note that this indexing does not add any labels to the final placement as the region of each rook can be uniquely determined.

The interior rooks have $k-j$ distinct values for each coordinate, all of which are different from~1. Thus, each coordinate has
\begin{equation*}
    q=n-(k-j)=m+j-1
\end{equation*}
available values, with exactly one value equal to~1. The available sets may vary between coordinates, but their cardinalities are all equal.

We assign the coordinates to the new $j$ rooks independently. For the coordinate $x_d$, all $j$ rooks are free and their values must be distinct. Thus, there are $\fall{q}{j}$ possibilities.

For each coordinate $a\in U$, exactly one specified rook has value~1. The remaining $j-1$ rooks take distinct values from the $q-1$ available non-1 values, and there are $\fall{q-1}{j-1}$ such assignments.

For each coordinate $a\in [d-1]\setminus U$, there is no rook with value~1. All $j$ rooks are assigned distinct values from the $q-1$ available non-1 values. There are $\fall{q-1}{j}$ such assignments.

Since every resulting cell belongs to its prescribed exclusive region, the board condition is satisfied automatically. The total number of extensions is therefore
\begin{equation}\label{eq:coordinate-count}
    \fall{q}{j} \left[\fall{q-1}{j-1}\right]^\sigma \left[\fall{q-1}{j}\right]^{d-1-\sigma}.
\end{equation}
Observe that
\begin{align*}
    \fall{q-1}{j-1} &= (q-1)(q-2)\cdots(q-j+1) \\
    & = (m+j-2)(m+j-3)\cdots(m) = \prod_{a=0}^{j-2}(m+a).
\end{align*}
Properties of the falling factoring yield the following identities:
\begin{equation*}
    % \fall{q}{j}=q(q-1)\cdots(q-j+1)=q \fall{q-1}{j-1}
    \fall{q}{j}=q \fall{q-1}{j-1}\quad\text{and}\quad 
    % \fall{q-1}{j}=(q-1)(q-2)\cdots(q-j)=(m+j-2)(m+j-3)\cdots(m-1)=(m-1)\fall{q-1}{j-1}
    \fall{q-1}{j}=(m-1)\fall{q-1}{j-1}.
\end{equation*}
Therefore, substituting these identities into~\eqref{eq:coordinate-count} gives us
\begin{equation*}
    q \left[\fall{q-1}{j-1}\right]^d (m-1)^{d-1-\sigma},
\end{equation*}
which is precisely~\eqref{eq:independence}, establishing the result.

\end{proof}

\begin{remark}
This proof also covers for the edge case $j=1$, when $\fall{q-1}{0}=1$ as well as $m=1$. For $m=1$, a boundary partition only contributes for $\sigma=d-1$, since otherwise one of the first $d-1$ coordinates would require $j$ distinct values from a set of size $j-1$.
\end{remark}

\begin{example}
For $d=4$, $j=2$, and $\sigma=3$, the three partition types are
\begin{equation*}
    \{\{1\},\{2,3\}\},\qquad \{\{2\},\{1,3\}\},\qquad \{\{3\},\{1,2\}\}.
\end{equation*}
Each contributes $(m+1)m^4$ to the overall count for a combined contribution of $3(m+1)m^4$.
\end{example}

\subsection{Proof of Theorem~\ref{thm:main}}
Finally, the above counting steps can now be combined to prove the uniform recurrence; we simply sum the extension count over all partition types and then over the possible numbers of boundary rooks.

\begin{proof}[Proof of Theorem~\ref{thm:main}]
Partition all placements of $k$ rooks according to the number of $j$ rooks on $\mathcal{B}_n^{(d)}$. By the pairwise disjointness of the $S_x$ sets associated with the boundary rooks, we have
\begin{equation*}
    0\leq j\leq\min(k,d-1).
\end{equation*}
If $j=0$, every rook lies in the translated interior. This contributes $r_k(\Gamma_{n-1}^{(d)})$.

Suppose $j\geq 1$. There are $r_{k-j}(\Gamma_{n-1}^{(d)})$ ways to choose the placements of the interior rooks. 
We can then count the number of extensions for each choice of the interior rooks by grouping the possible boundary partition types according to $\sigma$.
Then, using Equations~\eqref{eq:partition-count} and~\eqref{eq:independence}, we get the extension count:
\begin{equation*}
    \sum_{\sigma=j}^{d-1} \binom{d-1}{\sigma} S(\sigma,j)P_j^{(d)}(m,\sigma).
\end{equation*}
This is precisely $C_j^{(d)}(m)$.

Note that there is no factor of $j!$ because the blocks form an unordered partition, and the rook associated with each block is uniquely determined by its exclusive region. Thus, each placement is therefore counted exactly once. Summing over $j$ proves the recurrence.

\end{proof}
%----------------------------------------------------------------------
\section{Finite differences and coefficient identities}\label{sec:coefficients}
%----------------------------------------------------------------------
By utilizing the uniform recurrence in Theorem~\ref{thm:main}, we can extend the coefficient patterns from dimensions three through five to arbitrary dimension. Furthermore, we can describe the Stirling number factor as a sum of finite differences which allow us to develop general formulas that describe the degree, leading coefficient, and special values of the recurrence coefficients. These combinatorial interpretations help us see how the previously observed patterns from the lower-dimensional examples relate to each other.

To begin, we separate the formula in~\eqref{eq:main-coefficients} into two factors by introducing the notation
\begin{equation}\label{eq:structural-factor}
    F_j^{(d)}(m) = (m+j-1)\prod_{a=0}^{j-2}(m+a)^d\; \text{ and}
\end{equation}
\begin{equation}\label{eq:stirling-factor}
    B_j^{(d)}(m) = \sum_{\sigma=j}^{d-1} \binom{d-1}{\sigma} S(\sigma,j)(m-1)^{d-1-\sigma},
\end{equation}
so that $C_j^{(d)}=F_j^{(d)}B_j^{(d)}$ for $j\geq 1$.

Let $\Delta$ denote the forward-difference operator, i.e.,
\begin{equation*}
    (\Delta f)(x)=f(x+1)-f(x),
\end{equation*}
and for an integer $j\ge 0$, we write $\Delta^j$ for the operator obtained by
applying the $\Delta$ operator repeatedly $j$ times; that is,
\begin{equation*}
    \Delta^0 f = f, \qquad \Delta^{j} f = \Delta\bigl(\Delta^{j-1} f\bigr) \text{ for } j \geq 1.
\end{equation*}

% Equivalently, $(\Delta^j f)(x)=\sum_{i=0}^{j}(-1)^{j-i}\binom{j}{i}f(x+i)$.

\begin{proposition}[Finite-difference formula]\label{prop:finite-difference}
For $d\geq 3$ and $1\leq j\leq d-1$,
\begin{equation}\label{eq:finite-difference}
    B_j^{(d)}(m) = \frac{1}{j!} \sum_{t=0}^{j} (-1)^{j-t}\binom{j}{t}(m-1+t)^{d-1} = \frac{1}{j!} \left.\Delta^j x^{d-1}\right|_{x=m-1}.
\end{equation}
Consequently,
\begin{equation*}%\label{eq:coefficient-alternating}
    C_j^{(d)}(m) = \frac{m+j-1}{j!} \prod_{a=0}^{j-2}(m+a)^d\sum_{t=0}^{j}(-1)^{j-t}\binom{j}{t}(m-1+t)^{d-1}.
\end{equation*}
\end{proposition}

\begin{proof}
Use the well-known Stirling numbers inclusion-exclusion identity
\begin{equation*}
    S(\sigma,j) = \frac{1}{j!} \sum_{t=0}^{j} (-1)^{j-t}\binom{j}{t}t^\sigma.
\end{equation*}
Since $S(\sigma,j)=0$ for $0\leq\sigma<j$, the sum in~\eqref{eq:stirling-factor} can be extended to $\sigma=0$. 
Next, by interchanging the order of the two finite sums and utilizing the binomial theorem, we obtain
\begin{align*}
    B_j^{(d)}(m) &= \frac{1}{j!} \sum_{t=0}^{j} (-1)^{j-t}\binom{j}{t} \sum_{\sigma=0}^{d-1} \binom{d-1}{\sigma} t^\sigma(m-1)^{d-1-\sigma}\\
    &= \frac{1}{j!} \sum_{t=0}^{j} (-1)^{j-t}\binom{j}{t} (m-1+t)^{d-1}.
\end{align*}
Note that the second equality in~\eqref{eq:finite-difference} is just the standard expansion of the iterated forward difference.

\end{proof}

\begin{corollary}[Degree and leading coefficient]\label{cor:degree}
For $1\leq j\leq d-1$, the polynomial $C_j^{(d)}(m)$ has degree $j(d-1)$ with leading coefficient $\binom{d-1}{j}$. Moreover,
\begin{equation*}
    C_j^{(d)}(m)\in\mathbb{Z}_{\geq 0}[m-1].
\end{equation*}
\end{corollary}

\begin{proof}
The term $\sigma=j$ in~\eqref{eq:stirling-factor} has degree $d-1-j$ and leading coefficient $\binom{d-1}{j}$. All other terms have smaller degree. The structural factor has degree $1+d(j-1)$ with leading coefficient~1. Their product therefore has degree
\begin{equation*}
    (d-1-j)+1+d(j-1)=j(d-1)
\end{equation*}
and the stated leading coefficient.

To prove the last assertion, take $u=m-1$. Note that each factor in~\eqref{eq:structural-factor} is a power of a polynomial $u+c$ with $c$ a nonnegative integer. Also, all coefficients in~\eqref{eq:stirling-factor} are nonnegative integers for $u$.

\end{proof}

\begin{corollary}[Boundary coefficients]\label{cor:boundary-coefficients}
For $d\geq 3$,
\begin{align}
    \label{eq:C1}C_1^{(d)}(m) &= m\bigl(m^{d-1}-(m-1)^{d-1}\bigr),\\
    \label{eq:Clast} C_{d-1}^{(d)}(m) &= (m+d-2)\prod_{a=0}^{d-3}(m+a)^d.
\end{align}
More generally, for $1\leq j\leq d-1$,
\begin{equation}\label{eq:C-at-one}
    C_j^{(d)}(1) = j\bigl((j-1)!\bigr)^dS(d-1,j).
\end{equation}
\end{corollary}

\begin{proof}
Equation~\eqref{eq:C1} follows directly from Proposition~\ref{prop:finite-difference} with $j=1$. For $j=d-1$, the sum in~\eqref{eq:stirling-factor} consists only of the term $\sigma=d-1$, which equals~1. This proves~\eqref{eq:Clast}.

When $m=1$, all terms in~\eqref{eq:stirling-factor} vanish except $\sigma=d-1$. The structural factor becomes $j((j-1)!)^d$, proving~\eqref{eq:C-at-one}.

\end{proof}

\begin{remark}
Specializing Theorem~\ref{thm:main} to $k=n$ gives
\begin{equation}\label{eq:diagonal-from-triangle}
    r_n(\Gamma_n^{(d)}) = \sum_{j=1}^{\min(n,d-1)} j\bigl((j-1)!\bigr)^dS(d-1,j) r_{n-j}(\Gamma_{n-1}^{(d)}).
\end{equation}
Thus, the diagonal can be constructed using only a bounded number of entries in close proximity to the diagonal entry of the previous row; however, of course, this does not constitute a recurrence which involves only diagonal entries.
\end{remark}

%----------------------------------------------------------------------
\section{Low-dimensional specializations}\label{sec:examples}
%----------------------------------------------------------------------

The following specializations show how our uniform formula recovers previously established recurrences for dimensions three through five. To continue demonstrating the pattern, we also produce the six-dimensional case using the same approach. As before, we let $m=n-k+1$ throughout this section.

\subsection{Dimensions three through five}

For $d=3$, the nonconstant recurrence coefficients are
\begin{equation*}
    C_1^{(3)}(m)=m(2m-1), \text{ and } C_2^{(3)}(m)=(m+1)m^3.
\end{equation*}

Hence, Theorem~\ref{thm:main} gives
\begin{align*}
    r_k(\Gamma_n^{(3)}) &= r_k(\Gamma_{n-1}^{(3)}) \\
    & \quad + m(2m-1)r_{k-1}(\Gamma_{n-1}^{(3)}) \\
    & \quad + (m+1)m^3r_{k-2}(\Gamma_{n-1}^{(3)}),
\end{align*}
recovering the recurrence from~\cite{Alayont_Krzywonos_2016}.

For $d=4$, the coefficients are
\begin{align*}
    C_1^{(4)}(m)&=m(3m^2-3m+1),\\
    C_2^{(4)}(m)&=3(m+1)m^5,\\
    C_3^{(4)}(m)&=(m+2)(m+1)^4m^4.
\end{align*}
For example,
\begin{equation*}
    B_2^{(4)}(m)=3(m-1)+3=3m.    
\end{equation*}

For $d=5$, they are
\begin{align*}
    C_1^{(5)}(m)&=m\bigl(m^4-(m-1)^4\bigr),\\
    C_2^{(5)}(m)&=(m+1)m^5(6m^2+1),\\
    C_3^{(5)}(m)&=(m+2)(m+1)^5m^5(4m+2),\\
    C_4^{(5)}(m)&=(m+3)(m+2)^5(m+1)^5m^5.
\end{align*}
Here, the two intermediate sums simplify to
\begin{align*}
    B_2^{(5)}(m) &= 6(m-1)^2+12(m-1)+7 = 6m^2+1,\\
    B_3^{(5)}(m) &= 4(m-1)+6 = 4m+2.
\end{align*}
Note that the four- and five-dimensional formulas agree with the recurrences established in~\cite{Alayont_Loewen_Zadorozhnyy_2025}.

\subsection{Dimension Six}\label{sec:6d}
The following corollary establishes the recurrence for dimension six, extending the previously known low-dimensional cases.

\begin{corollary}\label{cor:six}
    For $n\geq 1$ and $0\leq k\leq n$,
    \begin{align*}
        r_k(\Gamma_n^{(6)}) &= r_k(\Gamma_{n-1}^{(6)}) \\
        &\quad+ m(5m^4-10m^3+10m^2-5m+1)r_{k-1}(\Gamma_{n-1}^{(6)}) \\
        &\quad+ 5(m+1)m^7(2m^2+1)r_{k-2}(\Gamma_{n-1}^{(6)}) \\
        &\quad+ 5(m+2)(m+1)^6m^6(2m^2+2m+1)r_{k-3}(\Gamma_{n-1}^{(6)}) \\
        &\quad+ 5(m+3)(m+2)^6(m+1)^7m^6r_{k-4}(\Gamma_{n-1}^{(6)}) \\
        &\quad+ (m+4)(m+3)^6(m+2)^6(m+1)^6m^6r_{k-5}(\Gamma_{n-1}^{(6)}),
    \end{align*}
    where $m=n-k+1$ and negative rook indices contribute zero.
\end{corollary}

\begin{proof}
The intermediate Stirling sums are
\begin{align*}
    B_2^{(6)}(m) &= 10(m-1)^3+30(m-1)^2+35(m-1)+15 \\
    &= 5m(2m^2+1),\\
    B_3^{(6)}(m) &= 10(m-1)^2+30(m-1)+25 \\
    &= 5(2m^2+2m+1),\\
    B_4^{(6)}(m) &= 5(m-1)+10 \\
    &= 5(m+1).
\end{align*}
Multiplying by the corresponding structural factors yields the coefficients of $r_{k-2}$, $r_{k-3}$, and $r_{k-4}$, and the remaining coefficients follow from Corollary~\ref{cor:boundary-coefficients}.

\end{proof}

The first few rows of the resulting rook-number triangle are recorded in Table~\ref{tab:six}. Observe that the first row accounts for the edge case $r_0(\Gamma_0^{(6)})=1$. The diagonal entries (highlighted in bold) align with OEIS sequence~\href{https://oeis.org/A065756}{A065756} from~\cite{oeis065756_gb_6d} and will be interpreted in Section~\ref{sec:gandhi}.

\begin{table}[ht]
    \centering
    \small
    \setlength{\tabcolsep}{4pt}
    \begin{tabular}{c c c c c c c}
        \toprule
        $n\backslash k$ & 0 & 1 & 2 & 3 & 4 & 5\\
        \midrule
        0 & \textbf{1} & & & & &\\
        1 & 1 & \textbf{1} & & & &\\
        2 & 1 & 63 & \textbf{31} & & &\\
        3 & 1 & 696 & 21217 & \textbf{6721} & &\\
        4 & 1 & 3820 & 1292845 & 25479615 & \textbf{5850271} &\\
        5 & 1 & 14325 & 26743325 & 5884739700 & 84597410401 & \textbf{15060446401}\\
        \bottomrule
    \end{tabular}
    \caption{Rook numbers $r_k(\Gamma_n^{(6)})$ for $0\leq k\leq n\leq 5$. Diagonal entries are in bold.}
    \label{tab:six}
\end{table}

As an additional elementary check on the first column of nonempty placements, the board definition gives
\begin{equation}\label{eq:single-rook}
    r_1(\Gamma_n^{(d)}) = n^d-\sum_{t=0}^{n-1}t^{d-1}.
\end{equation}
Indeed, for fixed $x_d=s$, an excluded cell has all of its first $d-1$ coordinates in $\{s+1,\ldots,n\}$. There are $(n-s)^{d-1}$ such cells.

Tables analogous to Table~\ref{tab:six} for dimensions $d=7$, $8$, and $9$ are given in Appendix~\ref{appendix_coefficients}.

%----------------------------------------------------------------------
\section{An injection model for generalized Gandhi polynomials}
\label{sec:gandhi}
%----------------------------------------------------------------------
The uniform recurrence determines the full rook-number triangle, but its diagonal specialization still involves off-diagonal entries, as seen in~\eqref{eq:diagonal-from-triangle}. To explain the diagonal sequence directly, we enlarge the permutation model to injections into larger sets. These additional codomain values supply the parameter for a bijective proof of the generalized Gandhi recurrence and hence of the known diagonal identification.

In order to make the dimensional dependence explicit, we define a family of polynomials $A_n^{(d)}(z)$ by
\begin{equation}\label{eq:gandhi-definition}
    A_0^{(d)}(z)=1, \text{ and } A_{n+1}^{(d)}(z) = z^{d-1}A_n^{(d)}(z+1) - (z-1)^{d-1}A_n^{(d)}(z), \quad n\geq 0.
\end{equation}
Thus,
\begin{equation*}
    A_1^{(d)}(z)=z^{d-1}-(z-1)^{d-1}.
\end{equation*}
These are the generalized Gandhi polynomials of degree $d-1$; see~\cite{Domaratzki2004CombinatorialIO} for more details. 
Note that our superscript records the board dimension instead of the exponent.

\subsection{Restricted injections}
We keep the inequality defining the board, using the last coordinate of a board element to identify the domain of an injection. The case where the injections are permutations gives us full rook placements on the board.

For $n\in\Nzero$ and an integer $r\geq 1$, let $\mathcal{A}_n^{(d)}(r)$ be the set of tuples $(\phi_1,\ldots,\phi_{d-1})$ such that
\begin{equation*}
    \phi_a:[n]\hookrightarrow[n+r-1]
\end{equation*}
is injective for each $a$, and
\begin{equation}\label{eq:injection-condition}
    \min_{1\leq a\leq d-1}\phi_a(i)\leq i \qquad\text{for every }i\in[n].
\end{equation}
There is one tuple of empty injections when $n=0$. Write
\begin{equation*}
    N_n^{(d)}(r)=|\mathcal{A}_n^{(d)}(r)|.
\end{equation*}

\begin{lemma}\label{lem:diagonal-injection}
    For every $d\geq 3$ and $n\in\Nzero$,
\begin{equation*}
    r_n(\Gamma_n^{(d)})=N_n^{(d)}(1).
\end{equation*}
\end{lemma}

\begin{proof}
For $r=1$, each injection is a permutation of $[n]$. A placement of $n$ nonattacking rooks on $\Gamma_n^{(d)}$ uses every value of every coordinate exactly once. In particular, for each $i\in[n]$, there is a unique rook with last coordinate $i$. Let $\phi_a(i)$ be its $a^{\text{th}}$ coordinate. The nonattacking condition makes each $\phi_a$ a permutation, and the board condition is exactly~\eqref{eq:injection-condition}.

Conversely, a tuple in $\mathcal{A}_n^{(d)}(1)$ gives the placement
\begin{equation*}
    \left\{(\phi_1(i),\ldots,\phi_{d-1}(i),i):i\in[n]\right\}.
\end{equation*}
These two operations are inverse to each other, establishing the result. Note this also holds for the edge case $n=0$: both sets consist of corresponding empty objects.

\end{proof}

\subsection{The deletion bijection}
To establish the Gandhi polynomials with the counts $N_n^{(d)}(r)$, we only have to show the recurrence~\eqref{eq:gandhi-definition} holds for the Gandhi polynomials. We do so by extending the domain by one element and by using a deletion bijection to count the extensions which violate the restriction.

\begin{theorem}\label{thm:injection-gandhi}
For every $d\geq 3$, $n\in\Nzero$, and integer $r\geq 1$,
\begin{equation}\label{eq:injection-gandhi}
    N_{n+1}^{(d)}(r) = r^{d-1}N_n^{(d)}(r+1) - (r-1)^{d-1}N_n^{(d)}(r).
\end{equation}
Consequently,
\begin{equation*}
    N_n^{(d)}(r)=A_n^{(d)}(r).
\end{equation*}
\end{theorem}

\begin{proof}
Fix
\begin{equation*}
    \psi=(\psi_1,\ldots,\psi_{d-1}) \in\mathcal{A}_n^{(d)}(r+1).
\end{equation*}
Each $\psi_a$ maps $[n]$ injectively into $[n+r]$, leaving exactly $r$ unused values. Choose one such value, say $b_a$, for each coordinate and extend $\psi_a$ by setting
\begin{equation*}
    \phi_a(n+1)=b_a.
\end{equation*}
There are $r^{d-1}$ extensions of each tuple $\psi$. Note that all conditions at positions $i\leq n$ remain satisfied.

An extension is not valid precisely when $b_a>n+1$ for every $a$. Thus, we construct a bijection between these invalid extensions and
\begin{equation}\label{eq:bad-product}
    \mathcal{A}_n^{(d)}(r)\times[r-1]^{d-1}.
\end{equation}

If $r=1$, both the set of invalid extensions and~\eqref{eq:bad-product} are empty. 

Now, assume that $r\geq 2$. Given an invalid extension, define $c_a=b_a-n-1\in[r-1]$. Since $b_a\notin\im(\psi_a)$, we can delete $b_a$ from the codomain and close the resulting gap:
\begin{equation*}%\label{eq:close-gap}
    \chi_a(i)=
    \begin{cases}
        \psi_a(i),&\psi_a(i)<b_a,\\
        \psi_a(i)-1,&\psi_a(i)>b_a.
    \end{cases}
\end{equation*}
Then, $\chi_a:[n]\hookrightarrow[n+r-1]$ is injective. Since $\chi_a(i)\leq\psi_a(i)$, it follows that~$\chi_a$ satisfies~\eqref{eq:injection-condition}. The forward map sends the invalid extension to $(\chi_a,(c_1,\ldots,c_{d-1}))$.

Conversely, take an element of~\eqref{eq:bad-product}. Set $b_a=n+1+c_a$ and open a gap at $b_a$:
\begin{equation}\label{eq:open-gap}
    \psi_a(i)=
    \begin{cases}
        \chi_a(i),&\chi_a(i)<b_a,\\
        \chi_a(i)+1,&\chi_a(i)\geq b_a.
    \end{cases}
\end{equation}
The resulting map is an injection into $[n+r]$ and does not use $b_a$.

It remains to check the restriction at each $i\leq n$. Choose $a$ such that $\chi_a(i)\leq i$. Since
\begin{equation*}
    \chi_a(i)\leq i\leq n<b_a,
\end{equation*}
this value is unchanged by~\eqref{eq:open-gap}, and, therefore, $\psi_a(i)\leq i$, as required. Appending $b_a$ at position $n+1$ now gives an invalid extension since every $b_a$ is greater than $n+1$.

The operations of closing and opening a gap are inverse to each other, so we can recover $b_a$ from $c_a$. This completes the bijection. The number of invalid extensions is given by
\begin{equation*}
    (r-1)^{d-1}N_n^{(d)}(r).
\end{equation*}
Subtracting this from the unrestricted extension count proves~\eqref{eq:injection-gandhi}.

Finally, the base case for the induction is $N_0^{(d)}(r)=A_0^{(d)}(r)=1$. The statement then follows by induction on $n$ using~\eqref{eq:gandhi-definition} and~\eqref{eq:injection-gandhi}.

\end{proof}
\begin{corollary}[Diagonal identification {\cite{Alayont2012ROOKNI}}]\label{cor:diagonal}
For every $d\geq 3$ and $n\in\Nzero$,
\begin{equation}\label{eq:diagonal-gandhi}
    r_n(\Gamma_n^{(d)})=A_n^{(d)}(1).
\end{equation}
\end{corollary}

\begin{proof}
This is an immediate consequence of combining Lemma~\ref{lem:diagonal-injection} with Theorem~\ref{thm:injection-gandhi}.

\end{proof}

\begin{remark}
In the notation of~\cite{Domaratzki2004CombinatorialIO}, the $p^{\text{th}}$ generalized Genocchi numbers are given $G_{2s}^{(p)}=A_{s-1}^{(p+1)}(1)$ in our notation. Thus,~\eqref{eq:diagonal-gandhi} can also be written as
\begin{equation*}
    r_n(\Gamma_n^{(d)})=G_{2n+2}^{(d-1)}.
\end{equation*}
The empty board accounts for the initial value at $n=0$. Keeping this index shift explicit avoids ambiguity between the board size and the index of a generalized Genocchi number.
\end{remark}

The Gandhi recurrence now provides a second way to compute the initial diagonal entries, allowing a direct comparison with the rook-number tables above.

\begin{corollary}[Initial diagonal values]\label{cor:initial-diagonal}
For every $d\geq 3$,
\begin{align*}
    r_1(\Gamma_1^{(d)}) &= 1,\\
    r_2(\Gamma_2^{(d)}) &= 2^{d-1}-1,\\
    r_3(\Gamma_3^{(d)}) &= 2^{d-1}\bigl(3^{d-1}-2^{d-1}\bigr) - \bigl(2^{d-1}-1\bigr).
\end{align*}
\end{corollary}

\begin{proof}
At $z=1$, the Gandhi recurrence gives
\begin{equation*}
    A_{n+1}^{(d)}(1)=A_n^{(d)}(2).
\end{equation*}
The first two values follow from $A_1^{(d)}(z)=z^{d-1}-(z-1)^{d-1}$. For the third,
\begin{align*}
    A_3^{(d)}(1) &= A_2^{(d)}(2) = 2^{d-1}A_1^{(d)}(3)-A_1^{(d)}(2),
\end{align*}
which gives the stated expression.

\end{proof}

For example, when $d=6$, $r_3(\Gamma_3^{(6)}) = 32(243-32)-31 = 6721$, as in Table~\ref{tab:six}.

%----------------------------------------------------------------------
\section{Concluding remarks}
\label{sec:conclusion}
%----------------------------------------------------------------------

Theorem~\ref{thm:main} gives a uniform recurrence for the whole rook-number triangle of a Genocchi board in any dimension. The coefficients in this recurrence separate the enumeration of the boundary partition types from the assignment of coordinates to the interior placement---the first part is given by Stirling numbers while the second part is independent of the geometry of the interior placement.

In order to understand better why the structural factor is the same for all types of partitions of given block size, a coordinatewise proof of Lemma~\ref{lem:independence} was introduced. The remaining factor can then be expressed in terms of a finite-difference formula.

The injection model provides a different explanation of the subtraction term in the generalized Gandhi recurrence. An invalid extension to a board can be described as a restricted injection into a smaller codomain, together with a choice of one deleted value in each of the coordinates. This gives an alternative explanation of the diagonal identification that was established using the complementary-board formula.

The aforementioned counting models, along with their formulas, raise some interesting questions for further study.
\begin{enumerate}[leftmargin=*]
    \item First, the recurrence for the full triangle can be used to study entries at a fixed distance from the diagonal, $r_{n-s}(\Gamma_n^{(d)})$, as $n$ varies. It would be of interest to compare such a recurrence to the polynomial families that arise from the injection model.
    \item Second, the proof of the uniform recurrence result relies on two compatible properties: translation preserves the interior board and fixing one of the first $d-1$ coordinates to~1 removes all remaining board inequalities on the boundary region corresponding to that coordinate. It would be interesting to identify other board families that have similar properties and use them to derive other recurrences. For boards defined using order statistics other than the minimum, the boundary regions will, in general, retain some inequalities, and thus the independence argument given above does not apply.
\end{enumerate}

Finally, the factorization of the coefficients and their nonnegative expansion in $m-1$ seem to suggest the possibility of establishing generating functions and inequalities for other entries among the full rook-number triangle instead of only for its diagonal. 

%----------------------------------------------------------------------

\clearpage
\bibliographystyle{plain} % We choose the "plain" reference style
\bibliography{bibliography} % Entries are in the refs.bib file

\makeatletter
\@setaddresses
\let\@setaddresses\relax % Prevents it from printing a second time at \end{document}
\makeatother
\clearpage

\begin{appendices}
\section{Rook-number tables in dimensions three through nine}\label{appendix_coefficients}

This appendix contains tables with the rook numbers $r_k(\Gamma_n^{(d)})$ (i.e., the coefficients of the rook polynomials of $\Gamma_n^{(d)}$) for $3\leq d\leq 9$ and $0\leq k\leq n\leq 5$. The tables are organized in the same way as in Table~\ref{tab:six}; that is, row $n$ consists of the rook numbers $r_0(\Gamma_n^{(d)}),\ldots,r_n(\Gamma_n^{(d)})$, where the entries with $k>n$ are 0 and thus left blank. The diagonal entries $r_n(\Gamma_n^{(d)})$ that count the placements of $n$ non-attacking rooks are bolded. The column $k=1$ is given in closed form by~\eqref{eq:single-rook}, providing a quick check on each table.

The tables in this appendix can be used to illustrate the full-triangle recurrence of Theorem~\ref{thm:main} as well as the identification of the diagonal entries in Corollary~\ref{cor:diagonal}. The diagonal entries for $d=3,4,5,6$ are listed in the OEIS sequences given below, and for $d=7,8,9$ are further specializations of the same generalized Genocchi family, as given in Appendix~\ref{appendix:7d9d}.

\subsection{3D}
Table~\ref{tab:coefficients_3d} lists the rook numbers of the three-dimensional boards $\Gamma_n^{(3)}$. As expected from previous discussions in~\cite{Alayont_Krzywonos_2016,Alayont_Loewen_Zadorozhnyy_2025}, the highlighted diagonal entries follow the sequence \href{https://oeis.org/A110501}{A110501}~\cite{oeis110501_gb_3d}: \emph{Unsigned Genocchi numbers (of first kind) of even index.} The diagonals in the higher-dimensional tables below generalize these numbers.

\begin{table}[ht]
    \centering
    \begin{tabular}{c|cccccc}
        $n \setminus k$ & 0 & 1 & 2 & 3 & 4 & 5 \\
        \hline
        0 & \bf{1} & & & & & \\
        1 & 1 & \bf{1} & & & & \\
        2 & 1 & 7      & \bf{3} & & & \\
        3 & 1 & 22     & 69     & \bf{17} & & \\    
        4 & 1 & 50     & 507    & 959     & \bf{155} &  \\
        5 & 1 & 95     & 2227   & 13964   & 18077    & \bf{2073}
    \end{tabular}
    \caption{Rook numbers $r_k(\Gamma_n^{(3)})$ for $0\leq k\leq n\leq 5$. Diagonal entries are in bold.}
    \label{tab:coefficients_3d}
\end{table}

\subsection{4D}
Table~\ref{tab:coefficients_4d} lists the rook numbers of $\Gamma_n^{(4)}$. Again as expected from~\cite{Alayont_Krzywonos_2016,Alayont_Loewen_Zadorozhnyy_2025}, the highlighted diagonal entries follow the sequence \href{https://oeis.org/A064624}{A064624}~\cite{oeis064624_gb_4d}: \emph{Generalization of the Genocchi numbers given by the Gandhi polynomials}
% \begin{equation*}
%     A(n+1,r) = r^3 \cdot A(n, r + 1) - (r - 1)^3 \cdot A(n, r); \; A(1,r) = r^3 - (r-1)^3.
% \end{equation*}
\begin{equation*}
    A_0^{(4)}(z) = 1, \text{ and } A_{n+1}^{(4)}(z) = z^3 A_n^{(4)}(z+1)
       - (z-1)^3 A_n^{(4)}(z), \quad n\geq 0,
\end{equation*}
% The terms of this sequence are the values $A(n,1)$, so $r_n(\Gamma_n^{(4)}) = A(n,1)$ for $n\geq 1$.
which is a specialized case of~\eqref{eq:gandhi-definition} with $d=4$. In particular, $A_1^{(4)}(z)=z^3-(z-1)^3$. By Corollary~\ref{cor:diagonal}, the diagonal entries satisfy $r_n(\Gamma_n^{(4)})=A_n^{(4)}(1)$ for $n\geq 0$.
\begin{table}[ht]
    \centering
    \begin{tabular}{c|cccccc}
        $n \setminus k$ & 0 & 1 & 2 & 3 & 4 & 5 \\
        \hline
        0 & \bf{1} & & & & & \\
        1 & 1 & \bf{1} & & & & \\
        2 & 1 & 15 & \bf{7} & & & \\
        3 & 1 & 72 & 505 & \bf{145} & & \\
        4 & 1 & 220 & 7525 & 33135 & \bf{6631} & \\
        5 & 1 & 525 & 55445 & 1207260 & 3778201 & \bf{566641}
    \end{tabular}
    \caption{Rook numbers $r_k(\Gamma_n^{(4)})$ for $0\leq k\leq n\leq 5$. Diagonal entries are in bold.}
    \label{tab:coefficients_4d}
\end{table}

\subsection{5D}

Table~\ref{tab:coefficients_5d} lists the rook numbers of $\Gamma_n^{(5)}$, and the pattern continues, again in agreement with~\cite{Alayont_Krzywonos_2016,Alayont_Loewen_Zadorozhnyy_2025}. The highlighted diagonal entries follow the sequence \href{https://oeis.org/A064625}{A064625}~\cite{oeis064625_gb_5d}: \emph{Generalization of the Genocchi numbers given by the Gandhi polynomials}

\begin{equation*}
    A_0^{(5)}(z) = 1, \text{ and } A_{n+1}^{(5)}(z)= z^4 A_n^{(5)}(z+1) - (z-1)^4 A_n^{(5)}(z), \quad n\geq 0.
\end{equation*}
% Again $r_n(\Gamma_n^{(5)}) = A(n,1)$ for $n\geq 1$; the only change from the four-dimensional case is that the exponent $3$ becomes $4$.
In particular, $A_1^{(5)}(z)=z^4-(z-1)^4$. Again, Corollary~\ref{cor:diagonal} gives $r_n(\Gamma_n^{(5)})=A_n^{(5)}(1)$ for $n\geq 0$; the exponent in the recurrence increases from $3$ to $4$.

\begin{table}[ht]
    \centering
    \begin{tabular}{c|cccccc}
        $n \setminus k$ & 0 & 1 & 2 & 3 & 4 & 5 \\
        \hline
        0 & \textbf{1} &  &  &  &  &  \\
        1 & 1 & \textbf{1} &  &  &  &  \\
        2 & 1 & 31 & \textbf{15} &  &  &  \\
        3 & 1 & 226 & 3345 & \textbf{1025} &  &  \\
        4 & 1 & 926 & 100875 & 954815 & \textbf{209135} &  \\
        5 & 1 & 2771 & 1245715 & 87547640 & 598789745 & \textbf{100482849}
    \end{tabular}
    % \caption{The coefficients of $r_k(\Gamma_n^{(5)}$).}
    \caption{Rook numbers $r_k(\Gamma_n^{(5)})$ for $0\leq k\leq n\leq 5$. Diagonal entries are in bold.}
    \label{tab:coefficients_5d}
\end{table}

\subsection{6D}
The six-dimensional case is treated in Section~\ref{sec:6d}, where the recurrence of Corollary~\ref{cor:six} produces the rook numbers in Table~\ref{tab:six}. As noted there, the diagonal entries follow the sequence \href{https://oeis.org/A065756}{A065756}~\cite{oeis065756_gb_6d}. They fit the same pattern: $r_n(\Gamma_n^{(6)}) = A(n,1)$, where $A$ is now given by the Gandhi recurrence with exponent $5$.

\subsection{Genocchi Numbers in Dimensions 7--9}\label{appendix:7d9d}

Corollary~\ref{cor:diagonal} also applies to the higher-dimensional cases, so the following tables extend the preceding comparisons beyond dimension six. For $d=7,8,9$, the highlighted diagonal entries in Tables~\ref{tab:coefficients_7d},~\ref{tab:coefficients_8d}, and~\ref{tab:coefficients_9d} are obtained by evaluating the generalized Gandhi polynomial recurrence
% \begin{equation*}
%     A(n+1,r) = r^{d-1} \cdot A(n, r+1) - (r-1)^{d-1} \cdot A(n, r), \; A(1,r) = r^{d-1} - (r-1)^{d-1},
% \end{equation*}
\begin{equation*}
    A_0^{(d)}(z) = 1, \text{ and } A_{n+1}^{(d)}(z) = z^{d-1} A_n^{(d)}(z+1) - (z-1)^{d-1} A_n^{(d)}(z), \quad n\geq 0.
\end{equation*}
% at $r = 1$, yielding $r_n(\Gamma_n^{(d)}) = A(n, 1)$ for $n\geq 1$. Unlike the cases $d\leq 6$, the resulting sequences appear to be new, and they have been submitted to the OEIS~\cite{oeis}.
In particular, $A_1^{(d)}(z)=z^{d-1}-(z-1)^{d-1}$, with exponents $6$, $7$, and $8$ in dimensions $7$, $8$, and $9$, respectively. By Corollary~\ref{cor:diagonal}, the highlighted diagonal entries in Tables~\ref{tab:coefficients_7d},~\ref{tab:coefficients_8d}, and~\ref{tab:coefficients_9d} are obtained by evaluating these polynomials at $z=1$:
\begin{equation*}
    r_n(\Gamma_n^{(d)})=A_n^{(d)}(1), \qquad n\geq 0.
\end{equation*}
The resulting sequences appear to be new, and they have been submitted to the OEIS~\cite{oeis}.

\begin{table}[ht]
    \centering
    \begin{tabular}{c|cccccc}
        $n \setminus k$ & 0 & 1 & 2 & 3 & 4 & 5 \\
        \hline
        0 & \textbf{1} &  &  &  &  &  \\
        1 & 1 & \bf{1} & & & & \\
        2 & 1 & 127 & \bf{63} & & & \\
        3 & 1 & 2122 & 131649 & \bf{42497} & & \\
        4 & 1 & 15590 & 16183587 & 653809919 & \bf{154324415} & \\
        5 & 1 & 73235 & 560465227 & 380225955704 & 11273511226817 & \bf{2078673245313}
    \end{tabular}
    % \caption{The coefficients of $r_k(\Gamma_n^{(7)}$).}
    \caption{Rook numbers $r_k(\Gamma_n^{(7)})$ for $0\leq k\leq n\leq 5$. Diagonal entries are in bold.}
    \label{tab:coefficients_7d}
\end{table}

% \subsection{8D}
\begin{table}[ht]
    \centering
    \resizebox{\textwidth}{!}{%
        \begin{tabular}{c|cccccc}
            $n \setminus k$ & 0 & 1 & 2 & 3 & 4 & 5 \\
            \hline
            0 & \textbf{1} &  &  &  &  &  \\
            1 & 1 & \bf{1} & & & & \\
            2 & 1 & 255 & \bf{127} & & & \\
            3 & 1 & 6432 & 806785 & \bf{263425} & & \\
            4 & 1 & 63220 & 199812085 & 16405117695 & \bf{3940253311} & \\
            5 & 1 & 371925 & 11579558405 & 24006609439140 & 1452789368918401 & \bf{273973317999361}
        \end{tabular}
    }
    % \caption{The coefficients of $r_k(\Gamma_n^{(8)}$).}
    \caption{Rook numbers $r_k(\Gamma_n^{(8)})$ for $0\leq k\leq n\leq 5$. Diagonal entries are in bold.}
    \label{tab:coefficients_8d}
\end{table}

% \subsection{9D}
\begin{table}[ht]
    \centering
    \resizebox{\textwidth}{!}{%
        \begin{tabular}{c|cccccc}
            $n \setminus k$ & 0 & 1 & 2 & 3 & 4 & 5 \\
            \hline
            0 & \textbf{1} &  &  &  &  &  \\
            1 & 1 & \textbf{1} &  &  &  &  \\
            2 & 1 & 511 & \textbf{255} &  &  &  \\
            3 & 1 & 19426 & 4907265 & \textbf{1613825} &  &  \\
            4 & 1 & 255326 & 2445757275 & 405980077055 & \textbf{98640535295} &  \\
            5 & 1 & 1880771 & 237077631715 & 1493832381997640 & 183359016043366145 & \textbf{35109020604482049} \\
        \end{tabular}
    }
    % \caption{The coefficients of $r_k(\Gamma_n^{(9)}$).}
    \caption{Rook numbers $r_k(\Gamma_n^{(9)})$ for $0\leq k\leq n\leq 5$. Diagonal entries are in bold.}
    \label{tab:coefficients_9d}
\end{table}
\end{appendices}

\hfill\hfill
\end{document}